\documentclass[preprint,11pt]{elsarticle}

\let\today\relax
\makeatletter
\def\ps@pprintTitle{%
    \let\@oddhead\@empty
    \let\@evenhead\@empty
    \def\@oddfoot{\footnotesize\itshape
         {} \hfill\today}%
    \let\@evenfoot\@oddfoot
    }
\makeatother

\usepackage{amssymb, amsmath}
\usepackage{kpfonts}
\usepackage{mathtools}
\newtheorem{teo}{Theorem}[section]

\newtheorem{lem}[teo]{Lemma}          %NumaralandÄ±rmayÄ± toerem ile birlikte yap, kÃ¶ÅŸeli paranteze numaralandÄ±rmanÄ±n yapÄ±lacaÄŸÄ± kÄ±sm yazÄ±yoruz.

 \usepackage{amsmath}
\begin{document}

\begin{frontmatter}

\title{A new lower bound for the number of conjugacy classes}

%\title{Bounding the number of conjugacy classes from below in terms of a prime }

\author[label1]{Burcu \c{C}{\i}narc{\i}}
\ead{bcinarci@txstate.edu}
\author[label1]{Thomas Michael Keller\corref{cor3} }
%\address [label1]{ Department of Mathematics, Texas State University, 601 University Drive, San Marcos, TX  78666}
%\address{Dedicated to the memory of Bertram Huppert (1927--2023)}
\address[label1]{ Department of Mathematics, Texas State University, 601 University Drive, San Marcos, TX \\ 
78666\\
Dedicated to the memory of Bertram Huppert (1927--2023)}
%\address [label2]{ Department of Marine Engineering, Piri Reis University,  Istanbul, Turkey 34940}

\ead{keller@txstate.edu} \cortext[cor3]{Corresponding author}

%\address[label3]{Department of Mathematics, Texas State University, 601 University Drive, San Marcos, TX 78666}

%\par\vspace*{.35\textheight}{\centering Dedicated to my parents\par}
%\dedicatory{Dedicated to the memory of N. N. (1930--2010)}

\begin{abstract}
In 2003, H\'{e}thelyi and
K\"{u}lshammer proposed that if $G$ is a finite group and $p$ is
a prime dividing the group order, then $k(G)\geq 2\sqrt{p-1}$, and they
proved this conjecture for solvable $G$ and showed that it is sharp for those
primes $p$ for which $\sqrt{p-1}$ is an integer. This initiated 
a flurry of activity, leading to many generalizations and variations of the result;
in particular, today the conjecture is known to be true for all finite groups. In this note,
we put forward a natural new and stronger conjecture, which is sharp for all primes $p$, and
we prove it for solvable groups, and when $p$ is large, also for arbitrary groups.
\end{abstract}


\end{frontmatter}

\section{Introduction}
\noindent The number of conjugacy classes of a finite group G is a 
fundamental parameter, usually denoted by $k(G)$, and has been studied intensely for decades. In this paper we are
once again concerned with bounding $k(G)$ from below.\\
\indent In 2003, H\'{e}thelyi and
K\"{u}lshammer \cite{hethelyi-kulshammer2000} proved that if $G$ is solvable and $p$ is
a prime dividing the group order, then $$k(G)\geq 2\sqrt{p-1}.$$
This bound is sharp, as can easily be seen as follows. Let $p$ be a
prime such that $p-1$ is a square. Let a cyclic group of order
$\sqrt{p-1}$ act frobeniusly on a cyclic group of order $p$, and let
$G$ be the corresponding semidirect product. Then $k(G) =
2\sqrt{p-1}$. It was conjectured that this result also holds for arbitrary finite groups, and indeed, through the combined efforts of several researchers in \cite{ke8}, \cite{he}, \cite{mor1} we have the following beautiful result.
\begin{teo}\label{x}
If $p$ is
a prime number and $G$ is a finite group of order
divisible by $p$, then $k(G)\geq 2\sqrt{p-1}$ with equality if and
only if $\sqrt{p-1}$ is an integer, $G = C_p \rtimes C_{\sqrt{p-1}}$
and $C_G(C_p) = C_p$.
\end{teo}
While this result is indeed very satisfying, it leaves some natural questions open.
It is not known whether there are infinitely many primes such that $\sqrt{p-1}$ is an integer, and this is a hard problem in number theory. (It is conjectured that there are
infinitely many such primes.) So, it is not known whether there are infinitely many groups, where the lower bound is attained. Also, one could still consider the bound sharp if
for other primes, we would find groups $G$ such that $k(G)$ equals the smallest integer
above $2\sqrt{p-1}$ (i.e., the ceiling of $2\sqrt{p-1}$), and then we do not know for which groups we have equality here.\\
\indent We therefore would like propose a new conjecture, which has a stronger conclusion than Theorem \ref{x}. In fact, it provides the exact best possible lower bound for every prime $p$
and also characterizes the groups for which the bound is sharp.\medskip\\
\noindent \textbf{Conjecture A.} 
\textsl{ Let $p$ be a prime. Let $a$ and $b$ be positive integers such that
$p-1 = ab$ and such that $|a-b|$ is minimal. Let $K= C_p \rtimes C_a$
such that $C_K(C_p) = C_p$,
and let $L=C_p \rtimes C_b$
such that $C_L(C_p) = C_p$.
Then
for any finite group $G$, which has order divisible by $p$ we have $k(G)\geq a+b$, with equality if and only if $G\cong K$
or $G\cong L$.}\medskip\\
\noindent A few remarks on Conjecture A are in order.\\
(1) First, observe that if $p\geq 5$, then $K$ and $L$ are Frobenius groups because
then $p-1$ will have a nontrivial factorization (i.e., both factors $a$ and $b$ are greater than 1). When $p=3$, then $\{a, b\}=\{1, 2\}$, and $K$ and $L$ are $C_3$ and $S_3$, respectively. When $p=2$, then $a=b=1$ and $K=L=C_2$. \\
(2) Also, clearly $a+b\geq 2\sqrt{p-1}$, with equality if and only if $\sqrt{p-1}$ is an integer, as would be expected in view of Theorem \ref{x}.\\
(3) Moreover, the groups $K$ and $L$ are non-isomorphic if and only if $a\ne b$, so that then we have two groups with minimal number of conjugacy classes. And we have $K\cong L$ if and only if $a=b=\sqrt{p-1}$ and we thus recover Theorem \ref{x} as a special case. \\
(4) We also note that for $p\geq 5$ the "worst" factorization of $p-1$ that is possible
clearly occurs when $(p-1)/2$ is a prime, in which case Conjecture A predicts that 
$k(G)\geq \frac{p-1}{2}+2=\frac{p+3}{2}$.
Note that this is significantly larger than what we get from Theorem \ref{x} in this case. Interestingly, the primes yielding this largest possible lower bound have also drawn significant interest in number theory. If $p$ is a prime such that $q=(p-1)/2$ is also a prime, then $q$ is called a {\it Sophie Germain prime}, and $p$ is called a {\it safe} prime. It is conjectured that there are infinitely many Sophie Germain primes.\\
(5) In general it is easy to check that if, for a given prime $p$, we let $a$ and $b$ be as in Conjecture A, then
\[2\sqrt{p-1}\leq a+b\leq \frac{p+3}{2}\]
(6) We also note another interesting feature of the new lower bound proposed by Conjecture A: It is not monotonic. To get the best lower bound using Theorem \ref{x}, one clearly has to identify the largest prime dividing the order of $G$. Not so with Conjecture A. If a group $G$ has order divisible by, say, 23 and 29, then Conjecture A, used with $p=23$, predicts that $k(G)\geq 13$, whereas used with $p=29$ it predicts that $k(G)\geq 11$. So, the smaller prime gives us a better bound (and Theorem \ref{x} with $p=29$ tells us that 
$k(G)\geq 11$ in this case).\\
(7) While the bound predicted by Conjecture A clearly is sharp for solvable groups, it can also be sharp for non-solvable groups. For example, it is well-known that 
\[k(\mbox{PSL}(2,p)) = \frac{p+5}{2} \] 
for any prime $p > 3$, and if $p$ is a safe prime as discussed above in (4), then Conjecture A predicts $k(\mbox{PSL}(2,p)) > a+b = (p+3)/2$, i.e. $k(\mbox{PSL}(2,p)) \geq (p+5)/2$. This shows that
the bound in Conjecture A for non-solvable groups is sharp for safe primes (of which there are expected to be infinitely many, as pointed out in (4)). It seems
likely that the groups $\mbox{PSL}(2,p)$ for safe $p$ are the only groups meeting the bound in the non-solvable case.\\
(8) In \cite{keller-moreto2023}, it is shown that there are (most likely) infinitely many primes such that
\[k(G) = \frac{p-1}{\log_2(p+1)}+2\log_2(p+1).\ \ \ (*)\]
This could potentially clash with Conjecture A since this value for $k(G)$ is smaller than $(p+3)/2$, and if $p$ were a safe prime as mentioned in (4), then we would have a contradiction. However, the examples realizing $(*)$ in
\cite{keller-moreto2023} are not counterexamples (they are even solvable, and so Theorem B below applies), and such examples seem fairly rare. If Conjecture A is true,
then this simply says that groups $G$ such that $k(G)$ is near $p/\log_2(p)$ for the largest prime divisor of $|G|$ can only occur for primes $p$ such that
$p-1$ has a divisor roughly of size $\log_2(p)$ and, in particular, cannot be a safe prime. This seems to be a curious insight.\\
\indent So far we have tried to ignite the reader's interest in Conjecture A. In the remainder of the paper, we will try to lend some credence to Conjecture A, by proving it in some important special cases as follows.\medskip\\
\noindent\textbf{Theorem B.} \textsl{Conjecture A is true for solvable groups.}\medskip\\
\noindent\textbf{Theorem C.} \textsl{Conjecture A is true for arbitrary finite groups and large $p$. More precisely, there exists a universal constant $C$ such that the following holds. Let $p>C$ and let $a$, $b$, $K$, $L$ be as in Conjecture A. If $G$ is a finite group of order divisible by $p$, then $k(G)\geq a+b$, with equality if and only if $G\cong K$ or $G\cong L$.}\medskip\\
\noindent In Section 2, we will prove and discuss these results, but we mention here that the proof of Theorem C 
depends on the recently proved McKay conjecture, and that
the constant $C$ is entirely determined by the case of
$p$-solvable groups. In other words, if one can prove 
Conjecture A for $p$-solvable groups and all primes $p$, then it will follow for all finite groups with the results in this
paper.\\
\indent Theorem \ref{x} has been strengthened in many ways by replacing $k(G)$ by smaller quantities, such as $|\textrm{Irr}_{p'}(G)|$,
the number of irreducible characters of $p'$ degree of $G$, or $k(B_0)$, the number of irreducible ordinary
characters in the principal $p$-block $B_0$ of $G$, see \cite{malle-maroti} and \cite{hung-schaeffer}, respectively. In the latter result, it is actually expected that the bound is true for any $p$-block $B$ instead of just the principal block.
 Even another
variation of this is found in \cite{hung-malle-maroti}.
Based on this we pose the following question.\medskip\\
\noindent\textbf{Problem D.} \textsl{For which of the variations of Theorem \ref{x}, where $k(G)$ is replaced by a smaller quantity, can the conclusion be replaced by the conclusion of Conjecture A? }\medskip\\
\noindent There is a good chance that the answer to the above question is yes for most of the quantities mentioned above. In fact, on our way to proving Theorem C we will need one
such result as asked for in Problem D, as follows.\medskip\\
\noindent\textbf{Theorem E.} \textsl{  There exists a universal constant $C$ such that the following holds. Let $p$ be a prime with $p>C$. Let $a$ and $b$ be positive integers such that
$p-1 = ab$ and such that $|a-b|$ is minimal. 
Then for any finite group $G$, which has order divisible by $p$ we have $|\textnormal{Irr}_{p'}(G)|\geq a+b$, where $|\textnormal{Irr}_{p'}(G)|$ 
is the number of irreducible complex characters of $G$ of degree not divisible by $p$.  }\medskip\\
The constant $C$ in Theorem E is the same as in Theorem C and comes
from proving Conjecture A for $p$-solvable groups; see Lemma \ref{p}
below.\\
Note that we currently cannot characterize the case of equality in Theorem E. In view of the results in 
\cite{malle-maroti}, characterizing equality in Theorem E might be a difficult task.
\section{Proofs of Theorems B, C, and E}
\noindent In this section,  we will prove Theorems B, C, and E. \medskip\\
\noindent\textbf{\textsl{Proof of Theorem B.}} Let $G$ be a counterexample of minimal order. Then as usual (see e.g., the proof of main theorem in \cite{hethelyi-kulshammer2000}) we have that $G=HV$ is a semidirect product of two subgroups $H$ and $V$ such that $V$ is a minimal normal subgroup of $p$-power order, and $H$ is a $p'$-group, which acts faithfully on $V$. First suppose that $p^2$ divides $|G|$. Then we invoke \cite{hethelyi-kulshammer2003} to conclude that \[k(G)\geq \frac{49p+1}{60}.\] 
If
$p\geq 5$, then it is easily checked that
$(49p+1)/60 > (p+3)/2$, which by the remark (5) following Conjecture A implies the desired conclusion, against $G$ being a counterexample. Since the cases $p=2$ and $p=3$ are trivial, we are done in case that $p^2$ divides $|G|$.\\
\indent So, now suppose that $p^2$ does not divide $|G|$. Then $G=HV$ is a Frobenius group with cyclic Frobenius kernel $V\cong C_p$ and cyclic Frobenius complement $H$, and if we write $|H|=c$ and $(p-1)/c=d$, then $c$ and $d$ are positive integers such that $cd=p-1$ and $k(G)=c+d$, as follows immediately from the well-known formula for the number of conjugacy classes of Frobenius groups. Hence by the definition of $K$ and $L$ it is clear that $G$ is not a counterexample. This contradiction concludes the proof.  $\hfill\Box$\medskip\\
\indent The proof of Theorem B is quite easy because for solvable groups we have the strong lower bound for $k(G)$ obtained in \cite{hethelyi-kulshammer2003}.
If for arbitrary groups of order divisible by $p^2$ we did know that 
$k(G)\geq \frac{1}{2}p+\frac{3}{2}$, then the proof of Conjecture A would work essentially the same as for solvable groups. However, for arbitrary groups of order divisible by $p^2$ we only have a lower bound of $cp$ by \cite{maroti-simion} for an unspecified constant $c>0$, which from the current proof certainly is quite small. Therefore for the moment we have to content ourselves to prove Conjecture A for large primes $p$ only and to appeal in our proof to the 
McKay conjecture, which, according to a recent announcement by B. Sp\"ath, now is
fully proved and thus henceforth will be referred to as {\it McKay theorem} in this paper. It might still be desirable to find a more elementary proof of Theorem C which does not depend on the McKay theorem. Before we can prove Theorem C, we need to establish a lemma and Theorem E. The following lemma proves Conjecture A for $p$-solvable groups.
\begin{lem}\label{p}
There exists a universal constant $C$ such that the following holds. Let $p>C$, and let $G$ be a finite $p$-solvable group. Let $a$, $b$, $K$, $L$ be as in Conjecture A. If $G$ has order divisible by $p$, then $k(G)\geq a+b$, with equality if and only if $G\cong K$ or $G\cong L$.
\end{lem}
\noindent\textbf{\textsl{Proof.}} Let $G$ be a counterexample of minimal order.
%By the McKay theorem, we have $|\textnormal{Irr}_{p'}(G)|=
%|\textnormal{Irr}_{p'}(N_G(P))|$ for a Sylow $p$-subgroup $P$ of $G$, so by minimality we may assume that $G=N_G(P)$,
%i.e., $P$ is normal in $G$. Then it is 
%well-known that $|\textnormal{Irr}_{p'}(G)|=|\textnormal{Irr}(G/P')|=k(G/P')$, and so by minimality we may further assume that $P'=1$ and $P$ is abelian.
%Note that at this point $|\textnormal{Irr}_{p'}(G)|=k(G)$, so we will write $k(G)$ instead of $|\textnormal{Irr}_{p'}(G)|$ from now on.\\
Now let $V$ be a minimal normal subgroup of $G$. If $p$ divides $|G/V|$, then $G/V$ is not a counterexample, and so we find that $k(G)\geq k(G/V)+1 > a+b$, and thus have a contradiction
to $G$ being a counterexample. Hence, $p$ does not divide $|G/V|$, and as $G$ is $p$-solvable, it follows that $V$ is an elementary abelian
$p$-group, which by Schur-Zassenhaus has a complement $H$ in $G$. Clearly $C_H(V)$ is
normal in $G$, and if $C_H(V)>1$, then $G/C_H(V)$ is not a counterexample, and
we again find that $k(G)\geq k(G/C_H(V))+1 > a+b$, a contradiction to $G$ being a
counterexample. Hence, $C_H(V)=1$. Thus, $G=HV$ is a semidirect product of two 
subgroups $H$ and $V$ such that $V$ is a minimal normal subgroup of $p$-power order, 
and $H$ is a $p'$-group, which acts faithfully on $V$.\\
\indent First suppose that $p^2$ does not divide $|G|$. Then $|V|=p$, and thus $G$ is solvable. Now we get a contradiction to $G$ being a minimal counterexample just as in the corresponding case in the 
proof of Theorem B. So, we now suppose that $p^2$ divides $|G|$. Then we are in the situation of 
\cite[Proposition 2.2]{maroti-simion}. Our situation is even more special in that the order of $H\cong G/V$ is not
divisible by $p$. We will thus follow the proof of \cite[Proposition 2.2]{maroti-simion}
closely and adjust it to our situation.\\
\indent First observe that it is well-known that $k(G) \geq k(H)+n(H,V)-1$,
where $n(H,V)$ denotes the number of orbits of $H$ on $V$. Let us assume that $p$ is sufficiently large. Every non-abelian (simple) composition factor of $H$ 
has order coprime to $p$ provided that it exists. There are the following possibilities for a 
non-abelian composition
factor $S$ of $H$: (i) $S$ is an alternating group; (ii) $S$ is a simple 
group of Lie type in characteristic different from $p$; 
(iii) $S$ is a sporadic simple group.
Observe that the case $S\cong \mbox{PSL}(2,p)$ in the proof of
\cite[Proposition 2.2]{maroti-simion} is not a case here because $p$
does not divide $|H|$. Suppose that a composition factor $S$ as in (i), (ii), or (iii) exists. Then we have 
$k(H)\geq k^*(S)$ by \cite[Lemma 2.5]{pyber1992}, where $k^*(S)$ denotes 
the number of Aut($S$)-orbits on the conjugacy classes of $S$.\\
\indent Let $S$ be an alternating group of degree $r\geq 5$. 
Since $|\mbox{Out}(S)|\leq 4$, we have $k^*(S)\geq k(S)/4$. Since 
$S$ is a normal subgroup of
index 2 in the symmetric group of degree $r$, we have $k(S)\geq \pi(r)/2$, where $\pi(r)$
denotes the number of partitions of $r$. We thus find that $k^*(S)\geq c_0^{\sqrt{r}}$
for some
constant $c_0 > 1$. If $r > (\log_2 p)^3$, then $k(G)\geq k(H)\geq k^*(S) > p$ for sufficiently large $p$, showing that $G$ is not a counterexample, contradiction.
Thus, we may
assume that if $S$ is an alternating composition factor of $H$, then $S$ has degree at most $(\log_2 p)^3$.\\
\indent Now let us write $|V|=p^n$ for an integer $n$. Then $H$ satisfies the hypothesis of \cite[Lemma 2.1]{maroti-simion}, and hence the group $H$ contains an abelian subgroup $A$ with
\[ |H:A|\leq (c_1 \log_2 p)^{7(n-1)}\leq (c_1 \log_2 p)^{7n},\ \ \ \ (*)\]
where $c_1>1$ is a universal constant.
%Since $p$ is sufficiently large, we have 
%$(c_1 \log_2 p)^{7}\leq p^{1/4}$ and thus %$|H:A|\leq |V|^{1/4}$.
Furthermore, by Ernest's result (see e.g. \cite[Exercise E17.3]{huppertcharacters}) we know that
$k(H)\geq k(A)/|H:A|=|A|/|H:A|$. Also clearly $n(H,V)\geq |V|/|H|$, and so we obtain
\[k(G)\geq k(H)+n(H,V)-1\geq \frac{|A|}{|H:A|}+\frac{|V|/|A|}{|H:A|}-1
=\frac{|A|+|V|/|A|}{|H:A|}-1,\]
and since the real function $g(x)=x+(|V|/x)$
takes its minimum in the interval
$[1, |V|]$ when $x=\sqrt{|V|}$ and by using 
$(*)$ we find that 
\[ k(G)\geq \frac{2 |V|^{\frac{1}{2}}}{|H:A|}-1 \geq
\frac{2 |V|^{\frac{1}{2}}}{(c_1 \log_2 p)^{7n}}-1 = 2\left(\frac{p^{\frac{1}{2}}}{(c_1 \log_2 p)^7}\right)^n -1 \]
Now note that since $p^2$ divides $|G|$, we know that 
$n\geq 2$. If even $n\geq 3$, then for sufficiently large $p$ we get
\[ k(G) \geq 2\left(\frac{p^{\frac{3}{2}}}{(c_1 \log_2 p)^{21}}\right)-1 >p\]
and thus $G$ is no counterexample, contradiction. This 
leaves us with the case $n=2$, i.e., $|V|=p^2$. By Theorem B
we may assume that $G$ (and thus $H$) is nonsolvable.
In this case, we have either $H/Z(H)\cong A_5$ or 
$H/Z(H)\cong S_5$ (given that $(|H|, |V|)=1$) by \cite[Section XII.260]{dickson} or
\cite[II, Hauptsatz 8.27]{huppert}.\\
\indent Now write $Z=Z(H)$ and consider $V_Z$, that is, $V$ viewed as a $Z$-module. If $V_Z$ is irreducible (i.e., $Z$ acts irreducibly on $V$), then
by \cite[II, Hilfssatz 3.11]{huppert} or \cite[Theorem 2.1]{manz} $G$ is solvable, a contradiction. Hence $V_Z$
is the direct sum of two $Z$-modules of order $p$, and since
clearly $Z$ acts Frobeniusly on $V$ (i.e., $ZV$ is a Frobenius group), this forces that $|Z|$ divides $p-1$.
Now clearly $k(H)\geq |Z|+1$ and 
\[n(H, V) \geq \frac{|V|}{|H|}\geq\frac{p^2}{120|Z|}=
\frac{p}{120}\ \frac{p}{|Z|} > \frac{p-1}{|Z|},\]
where the last inequality follows as we certainly may assume
that $p>120$.
Therefore we see that
\[ k(G)\geq k(H)+n(H,V)-1 > |Z|+ \frac{p-1}{|Z|}=c+d, \]
where $c=|Z|$ and $d=\frac{p-1}{|Z|}$. But since 
$cd=p-1$, from the definition of $a$ and $b$ we 
infer that $c+d\geq a+b$. Thus altogether 
$k(G) > a+b$, contradicting once again $G$
being a counterexample. This final contradiction 
completes the proof of the lemma.$\hfill\Box$\medskip\\
It seems that if desired, by careful inspection of the proof of Lemma \ref{p} (and Lemma 2.1 in \cite{maroti-simion}) a specific value of $C$ in Lemma \ref{p} could be determined.\\
\indent Next we prove Theorem E. This is the place where the McKay theorem comes into play and establishes the $a+b$-bound of Conjecture A, but not the characterization of equality, for a smaller quantity than $k(G)$, namely
$|\textnormal{Irr}_{p'}(G)|$.\medskip\\
\noindent\textbf{\textsl{Proof of Theorem E.}} Let $C$ be as in Lemma \ref{p} and let $p>C$.
Let $G$ be a counterexample of minimal order. By the McKay theorem, we have $|\textnormal{Irr}_{p'}(G)|=
|\textnormal{Irr}_{p'}(N_G(P))|$ for a Sylow $p$-subgroup $P$ of $G$, so by minimality we may assume that $G=N_G(P)$,
i.e., $P$ is normal in $G$. 
%Then it is 
%well-known that $|\textnormal{Irr}_{p'}(G)|=|\textnormal{Irr}(G/P')|=k(G/P')$ since $P/P'$ is an abelian and normal subgroup of $G/P'$. Hence by minimality we may further assume that $P'=1$, and so $P$ is abelian.
%Note that at this point $|\textnormal{Irr}_{p'}(G)|=|\textnormal{Irr}(G)|=k(G)$, so we will write $k(G)$ instead of $|\textnormal{Irr}_{p'}(G)|$ from now on. Assume for a moment
%that $O_{p'}(G)>1$. Since the order of $G/O_{p'}(G)$ is divisible by $p$ and $k(G)\geq k(G/O_{p'}(G))\geq |\textnormal{Irr}_{p'}(G/O_{p'}(G))|$, we get a contradiction by induction. Then we conclude that $O_{p'}(G)=1$, and hence $C_G(P)=P$ since $C_G(P)$ is a normal subgroup of $G$ and $O_{p'}(G)=1$. It follows from Schur-Zassenhaus theorem that $P$ has a complement $H$ in $G$, and $H$ acts faithfully on $P$. 
Thus, $G$ is a $p$-solvable group, and we get a contradiction $|\textnormal{Irr}_{p'}(G)|=k(G)\geq a+b$ by Lemma \ref{p}, which proves Theorem E. $\hfill\Box$\medskip\\
%It is even true if we replace $k(G)$ by $|\textnormal{Irr}_{p'}(G)|$  in the statement, where $|\textnormal{Irr}_{p'}(G)|$ 
%is the number of irreducible complex characters of $G$ of degree not divisible by $p$.
We finally can prove our main result.\medskip\\
\noindent\textbf{\textsl{Proof of Theorem C.}}  Let $C$ be as in Lemma \ref{p} and let $p>C$. By Theorem E, we have that  $k(G)\geq|\textrm{Irr}_{p'}(G)|\geq a+b$. Thus, we will prove the second part of the theorem. Let us assume that $k(G)=a+b$ and that $G$ is a counterexample. Since $a+b=k(G)=|\textrm{Irr}(G)|\geq |\textrm{Irr}_{p'}(G)|
%=|\textrm{Irr}_{p'}(N_G(P))|
\geq a+b$ by Theorem E, we have equality everywhere and hence $ |\textrm{Irr}_{p'}(G)|=|\textrm{Irr}(G)|$. Thus, it is well known that $P$ is a normal and abelian subgroup of $G$. 
In particular, $G$ is $p$-solvable, and we are done by Lemma \ref{p}.
%Then, since $p$ divides the order of $G/O_{p'}(G)$, we have $a+b=k(G)=|\textrm{Irr}(G)|\geq |\textrm{Irr}_{p'}(G/O_{p'}(G))|\geq a+b$ by Theorem E. Thus, we have equality everywhere, and  so $O_{p'}(G)=1$. By considering the fact that $C_G(P)=P$, we conclude that $G=HP$ and $H$ acts faithfully on $P$ for a complement $H$ of $P$. So $G$ is $p$-solvable, and we are done by Lemma \ref{p}. 
$\hfill\Box$ 
\section*{Acknowledgement} The authors would like to thank Attila Mar\'{o}ti and Iulian Simion for some fruitful discussions on the topics of this paper.
%This work   was done while the first author visited the second author as a Research Fellow, supported
 % by the Scientific and Technological Research Council of Turkey, at Texas State University.  She  would like to thank
 %the Department of Mathematics at Texas State University for its
 %hospitality, and T\"UB\.ITAK for granting her the research fellowship.\\

 %{\bf Conflict of interest statement:} On behalf of all authors, the corresponding author states that there is no conflict of interest. \\

 %{\bf Data availability statement:} Data sharing not applicable to this article as no datasets were generated or analysed during the current study.\\

\end{document}